\documentclass{amsart}

\newtheorem{theorem}{Theorem}
\newtheorem{lemma}[theorem]{Lemma}

\newcommand{\fix}{\operatorname{F}}
\newcommand{\least}{\operatorname{L}}

\begin{document}

\title[Group automorphisms with few periodic points]{Group automorphisms
with few\\and with many periodic points}

\author{Thomas Ward}
\address{School of Mathematics, University of East
  Anglia, Norwich NR4 7TJ, United Kingdom}
\email{t.ward@uea.ac.uk}


\date{June 16, 2003}



\begin{abstract}
For any $C\in[0,\infty]$ a compact
group automorphism $T:X\to X$ is constructed
with the property that
$$
\frac{1}{n}\log\vert\{x\in X\mid
T^n(x)=x\}\vert\longrightarrow C.
$$
This may be interpreted as a combinatorial
analogue of the (still open) problem of
whether compact group automorphisms exist
with any given topological entropy.
\end{abstract}

\maketitle

\section{Introduction}

One of the outstanding open problems in the dynamical properties
of compact group automorphisms concerns the existence of group
automorphisms with small entropy. Taking the infimum over all
compact group automorphisms $T$, and writing $h_{\rm top}(T)$ for
the topological entropy of $T$, is
\begin{equation}\label{infofentropy}
 \inf\{h_{\rm top}(T)\mid h_{\rm top}(T)>0\}>0?
\end{equation}
As pointed out in~\cite{MR49:10856}, the infimum being zero is
equivalent to the statement that for any $H\in(0,\infty]$ there is
an (ergodic) compact group automorphism with topological entropy
$H$. The infimum being positive implies that the set of possible
values of topological entropies of compact group automorphisms is
countable.

The problem in \eqref{infofentropy} turns out to be exactly
equivalent to Lehmer's Problem (see~\cite{MR49:10856}).
In~\cite{lehmer-1933}, Lehmer associated to any monic polynomial
$$
f(x)=\prod_{i=1}^{d}(x-\alpha_i)\in\mathbb Z[x]
$$
the integer sequence defined by
$\Delta_n(f)=\prod_{i=1}^{d}\vert\alpha_i^n-1\vert$, and asked how
fast this sequence grows. To avoid degeneracies, we assume that
$\Delta_n(f)\neq0$ for all $n\ge1$ (equivalently, $f$ does not
vanish on any root of unity). It turns out that the measure of
growth he used does not always exist, but a consequence of Baker's
Theorem is that the logarithmic growth rate
\begin{equation}\label{growthratetoral}
  m(f)=\lim_{n\to\infty}\frac{1}{n}\log\Delta_n(f)=\sum_{i=1}^{d}
\max\{\log\vert\alpha_i\vert,0\}
\end{equation}
exists
(see~\cite[Chap.~1]{MR2000e:11087} for the details). Lehmer's
problem is the following question: is
\begin{equation}\label{lehmerproblem}
  \inf\{m(f)\mid f\in\mathbb Z[x], m(f)>0\}>0?
\end{equation}
The equivalence of~\eqref{infofentropy} and~\eqref{lehmerproblem}
is proved in~\cite{MR49:10856}; to see why it is plausible notice
that if $T_f$ is the endomorphism of the $d$-torus $\mathbb T^d$
associated to the companion matrix of the polynomial $f$, then
$$h_{\rm top}(T_f)=m(f).$$
More is true. The sequence $\left(\Delta_n(f)\right)$ counts the
periodic points under the automorphism $T_f$,
\begin{equation}\label{toralperiodicpoints}
  \fix_n(T_f)=\left\vert\{x\in\mathbb T^d\mid
  T_f^n(x)=x\}\right\vert=\Delta_n(f)
\end{equation}
whenever either side is finite (ergodicity of $T_f$ is equivalent
to $\fix_n(T_f)$ being finite for all $n\ge1$).
Thus~\eqref{growthratetoral} means that the logarithmic growth
rate of the number of periodic points for an ergodic toral
automorphism coincides with the topological entropy. Indeed, any
continuous map on a compact metric space with sufficiently strong
specification properties will have this property
(see~\cite{MR96c:58055} for a detailed exposition). For compact
group automorphisms, expansiveness guarantees that the growth rate
of periodic points is the topological entropy (this is shown
in~\cite{MR92j:22013} in greater generality).

Without the assumption of expansiveness or ergodicity, we show
that in contrast to the conjectured answer `yes' to the
questions~\eqref{infofentropy} and~\eqref{lehmerproblem} above,
all possible logarithmic growth rates of periodic points arise for
compact group automorphisms.

\section{Periodic orbits}

The examples will be constructed by controlling the growth in the
number of orbits of length $n$, so the first step is to relate
this to the number of points of period $n$. Given any map $T:X\to
X$, define
\begin{equation}\label{defoffixn}
  \fix_n(T)=\left\vert\{x\in X\mid T^n(x)=x\}\right\vert
\end{equation}
to be the number of points of period $n$ and
\begin{equation}\label{defofleastn}
  \least_n(T)=\left\vert\{x\in X\mid\{T^j(x)\}_{j\in\mathbb
  N}=n\}\right\vert
\end{equation}
to be the number of points with least period $n$. The sequences
$\left(\fix_n(T)\right)$ and $\left(\least_n(T)\right)$ determine
each other via the relation~\eqref{basicsumformula} and its
M\"obius inversion (see~\cite[Sect.~11.2]{epsw}
or~\cite{MR2002i:11026}).

\begin{lemma}\label{fixleastrelation} For $C\in(0,\infty)$,
$$\lim_{n\to\infty}\frac{1}{n}\log\fix_n(T)=C$$
if and only if
$$\lim_{n\to\infty}\frac{1}{n}\log\least_n(T)=C.$$
\end{lemma}

\begin{proof}
This is shown in~\cite[Sect.~4]{MR2002i:11026}; a short proof is
included here. First notice that
\begin{equation}\label{basicsumformula}
  \fix_n(T)=\sum_{d\vert n}\least_d(T).
\end{equation}
If $\frac{1}{n}\log\least_n(T)\to C>0$ then, for $n$ large enough to
have $\fix_n(T)>0$,
\begin{eqnarray*}
\frac{1}{n}\log\least_n(T)\le\frac{1}{n}\log\fix_n(T)&=&
\frac{1}{n}\log\left(\sum_{d\vert n}\least_d(T)\right)\\
&\le&\frac{1}{n}\log n+\frac{1}{n}\log\max_{d\vert
n}\{\least_d(T)\}.
\end{eqnarray*}
For each such $n$, choose $\tilde{n}\in\{d\mid d\vert n,
\least_d(T)\ge\least_{d'}(T){\ }\forall{\ }d'\vert n\}$ so that
$\least_{\tilde{n}}(T)=\max_{d\vert n}\{\least_d(T)\}$ and
$\frac{\tilde{n}}{n}\le 1$. Now $\tilde{n}\to\infty$
as $n\to\infty$ since $C>0$. It follows that
\begin{eqnarray*}
\frac{1}{n}\log\least_n(T)\le\frac{1}{n}\log\fix_n(T)&\le& \frac{1}{n}\log
n+\frac{\tilde{n}}{n}\cdot
\frac{1}{\tilde{n}}\log\least_{\tilde{n}}(T)\\
&\le& \frac{1}{n}\log
n+\frac{1}{\tilde{n}}\log\least_{\tilde{n}}(T) \to C,
\end{eqnarray*}
so $\frac{1}{n}\log\fix_n(T)\to C$.

Now assume that $\frac{1}{n}\log\fix_n(T)\to C$. For $r\ge 1$,
\begin{equation*}
\fix_r(T)\ge\least_r(T)=-\sum_{d\vert r,d\neq
r}\least_d(T)+\fix_r(T)\ge\fix_r(T)-\sum_{d\vert r,d\neq
r}\fix_d(T).
\end{equation*}
Let $R$ be an upper bound for $\{\frac{1}{n}\log\fix_n(T)\mid
\fix_n(T)\neq0\}$; pick $\epsilon\in(0,3C)$. Choose $N$ so that
$$
r>N\implies e^{r(C-\epsilon)}\le\fix_r(T)\le e^{r(C+\epsilon)}.
$$
Then for $r>2N$,
\begin{eqnarray*}
\fix_r(T)\ge\least_r(T)&\ge&\fix_r(T)-\sum_{n=1}^{N}\fix_n(T)-
\sum_{n=N+1}^{\lfloor r/2\rfloor}\fix_n(T)\\
&\ge&
\fix_r(T)-\left(Ne^{NR}+(r/2-N)e^{r(C+\epsilon)/2}\right)\\
&\ge&
\fix_r(T)\left(1-Ne^{NR-r(C-\epsilon)}-(r/2-N)e^{-r(C-3\epsilon)/2}\right),
\end{eqnarray*}
and the bracketed expression converges to $1$ as $r\to\infty$.
Taking logs and dividing by $r$ gives the result.
\end{proof}

Write $\mathbb F_p=\mathbb Z/p\mathbb Z$ for the finite field with
$p$ elements. For any prime $p$, the field $\mathbb F_p$ has a
primitive root $g$ that generates the cyclic multiplicative group
$\mathbb F_p^{\ast}$.

\begin{theorem}\label{ratetheorem}
For any $C\in[0,\infty]$, there is a compact group automorphism
$T$ with $\fix_n(T)<\infty$ for all $n$ and with
$$
\frac{1}{n}\log\fix_n(T)\longrightarrow C.
$$
\end{theorem}

\begin{proof} If $C=0$ then we may take $X$ to be any finite group
(more interesting ergodic examples with $X=\widehat{\mathbb Q}$ or
$\widehat{\mathbb F_p(t)}$ may be found in~\cite{MR99b:11089} and
\cite{MR90a:28031}).

Assume now that $C=\infty$. Choose primes $p_1,p_2,\dots$ with the
following properties:
\begin{enumerate}
  \item $\frac{1}{n}\log p_n\longrightarrow\infty$, and
  \item $n\vert p_n-1$ for all $n\ge1$.
\end{enumerate}
This can be done since there are infinitely many primes in the
arithmetic progression $(1+nk)_{k\in\mathbb N}$ by Dirichlet's
theorem. For each $i$, let $g_i$ be a primitive root modulo $p_i$
and define a map $T_i:\mathbb F_{p_i}\to\mathbb F_{p_i}$ by
$T_i(x)=g_i^{(p_i-1)/i}x$ modulo $p_i$. Then $\fix_n(T_i)=p_i$ if
and only if $i\vert n$; in all other cases $\fix_n(T_i)=1$. Define
$X=\prod_{i=1}^{\infty}\mathbb F_{p_i}$ and the automorphism
$T=\prod_{i=1}^{\infty}T_i$. Then $\fix_n(T)\ge p_n$, so
$\frac{1}{n}\log\fix_n(T)\to\infty$. Notice that $\fix_n(T)$ is
finite for all $n$ since $p_n\to\infty$.

The argument above uses the fact that arithmetic progressions of
the form $(1+nk)_{k\in\mathbb N}$ contain infinitely many primes.
For $C<\infty$, a related but deeper fact is used.
Linnik~\cite{MR6:260b} has shown that the least prime congruent to
$1$ modulo $n$ is no larger than $n^A$ for some absolute constant
$A$. Later developments have culminated in Heath--Brown's
proof~\cite{MR93a:11075} that the bound may be replaced by
$Bn^{5.5}$ for an explicitly computable constant $B$.

The case $C\in(0,\infty)$ remains. By
Lemma~\ref{fixleastrelation}, it is enough to construct a group
automorphism $T$ with $\least_n(T)$ close to $e^{nC}$. By Linnik's
Theorem, choose primes $p_1,p_2,\dots$ with the following
properties:
\begin{enumerate}
  \item $n\vert p_n-1$ for all $n\ge1$, and
  \item $p_n\le n^A$.
\end{enumerate}
Let $g_n$ be a primitive root modulo $p_n$, and define
\begin{equation}\label{defineK}
  K_n=\left\lfloor\frac{nC}{\log p_n}\right\rfloor\mbox{ for }n\ge1,
\end{equation}
where $\lfloor x\rfloor$ is the largest integer less than or equal
to $x$. Notice that the second property of $p_n$ ensures that
$\frac{n}{\log p_n}\to\infty$ as $n\to\infty$. For each $n$,
let
$$
T_n:\left(\mathbb F_{p_n}\right)^{K_n}\to\left(\mathbb
F_{p_n}\right)^{K_n}
$$
be defined by
$$T_n(x_1,\dots,x_{K_n})=(g_n^{(p_n-1)/n}x_1,\dots,
g^{(p_n-1)/n}_nx_{K_n}).$$ Define the
compact group by
$$
X=\prod_{i=1}^{\infty}\left(\mathbb F_{p_i}\right)^{K_i},
$$
and the automorphism by $T=\prod_{i=1}^{\infty}T_i$.

The map $T_n^j$ is the identity if and only if $n\vert j$.
So $T_n^j(x_1,\dots,x_{K_n})=(x_1,\dots,x_{K_n})$ if and only if
$x_1=\dots= x_{K_n}=0$ or $n\vert j$.
It
follows that the number of points of least period $n$ under $T$ is
given by
\begin{equation}\label{leastperiodinT}
  \least_n(T)=p_n^{K_n}-1.
\end{equation}
Thus
\begin{eqnarray*}
\frac{1}{n}\log\least_n(T)&=&\frac{1}{n}\log(p_n^{K_n}-1)\\
&=&\frac{K_n}{n}\log(p_n^{K_n}-1)^{1/K_n}\\
&=&\frac{1}{n}\left\lfloor\frac{nC}{\log p_n}\right\rfloor
\log(p_n^{K_n}-1)^{1/K_n}
\longrightarrow C
\end{eqnarray*}
as $n\to\infty$ since $\frac{n}{\log p_n}\to\infty$. It follows
by Lemma~\ref{fixleastrelation} that
$$
\frac{1}{n}\log\fix_n(T)\longrightarrow C\mbox{ as }n\to\infty.
$$
\end{proof}

\section{Remarks}
\noindent(1) The only requirement in the proof is that
$K_n\to\infty$, so essentially the same argument allows the
construction of group automorphisms with any super-polynomial
growth rate, $\frac{\log\fix_n(T)}{\log n}\longrightarrow\infty$.

\noindent(2) A similar infinite product construction is used
in~\cite{MR1938222} to exhibit a group automorphism whose periodic
points count the Bernoulli denominators.

\noindent(3) The proof of Theorem~\ref{ratetheorem} constructs a
non-ergodic automorphism of a totally disconnected group. Does the
same result hold for ergodic automorphisms of connected groups?
On connected groups, the type of behaviour
seen in automorphisms with few periodic points seems to be
$\limsup_{n\to\infty}\frac{1}{n}\log\fix_n(T)=h_{\rm top}(T)$
and $\liminf_{n\to\infty}\frac{1}{n}\log\fix_n(T)=0$
(see~\cite{MR98k:22028}, \cite{MR99k:58152}).

\noindent(4) The dynamical zeta function associated to a map
$T$ with $\fix_n(T)<\infty$ for all $n\ge1$ is the formal
power series $\exp\sum_{n=1}^{\infty}
\frac{z^n}{n}\fix_n(T).$ Explicit examples of group
automorphisms with irrational zeta functions
are easy to find (see~\cite{MR99b:11089}; it should
be noted that the examples of this kind in~\cite{MR46:6400}
do not seem to be correct). For values of $C$ that are not
reciprocals of algebraic integers, the examples constructed
here cannot have rational zeta functions for a purely
arithmetic reason (the poles and zeros of a rational zeta
function can only occur at reciprocals of algebraic
integers by~\cite[Appendix]{MR42:6284} or~\cite[Part~VIII, Chap.~4,
No.~230]{MR57:5529}).


\end{document}